\documentclass[12pt]{amsart}

\def\withfig{0}
 
\ifnum\withfig=1
\usepackage{psfig}
\fi

\allowdisplaybreaks[4]

\begin{document}

\def\nsubseteq{{\not\subset}}
\def\vpp{{\vskip1cm}}
\def\vpd{{\vskip.9cm}}
\def\vnl{{\vskip.3cm}}
\def\vvv{{\vskip.6cm}}
\def\hhh{{\noindent \hskip1cm}}
\def\proof{{{\it Proof.}  ---\hskip.3cm}}
\def\proofend{{ \hskip.5cm QED}}
\def\da{{\downarrow}}
\def\OO{{\mathcal O}}
\def\FF{{\mathcal F}}
\def\EE{{\mathcal E}}
\def\PP{{\mathcal P}}
\def\PPP{{\mathbb P}}
\def\ZZ{{\mathbb Z}}
\def\QQ{{\mathbb Q}}
\def\CC{{\mathbb C}}
\def\II{{\mathcal I}}
\def\LL{{\mathcal L }}
\def\lra{\longrightarrow}
\def\ot{{\otimes}}
\def\ss{{\mathcal S}}
\def\sh{{\hat {\mathcal S}}}
\def\Pic{{\rm {Pic}^0}}
\def\dim{{\rm {dim}}}
\def\codim{{\rm {codim}}}
\def\alb{{\rm {alb}}}
\def \Cal{\mathcal} 
\def \Hom{{\mathrm{Hom}}}
\def \Supp{{\mathrm{Supp}}}
\def \closure#1{\overline{#1}}
\def \isom{\cong}
\def \embed{\hookrightarrow}

\begin{center} {\bf \Large Fourier transforms, generic vanishing 
theorems

\vskip.2cm
and polarizations of abelian varieties} 

\vskip.4cm
{By {\it Christopher D. Hacon} at {\it Utah}}

\vskip.4cm
-----------------------------
\end{center}
\footnotetext{{\it Key words and phrases:} Abelian variety, polarization, singularity}
\footnotetext{{\it 1991 Mathematics Subject Classification :} Primary 14K25; Secondary 14E05.}
\vskip.4cm
\noindent \hskip1cm {\bf Abstract.}\hskip.3cm The purpose of this paper is 
to give two applications of Fourier
transforms and generic vanishing theorems: 

--- we give a 
cohomological characterization of principal polarizations 

--- we prove that 
if $X$ an abelian variety and $\Theta $ a 
polarization of type $(1,...,1,2)$, then a general pair $(X,\Theta )$
is log canonical. 

\vvv
\begin{center} {\bf Introduction}
\end{center}
\vvv
\noindent \hskip1cm There is a well known connection between the geometry of
principally polarized abelian varieties (PPAVs), 
and the singularities of their theta divisors.
This was first discovered by Andreotti and Mayer, \cite {AM}
in their work on the Schottky locus, and has since found applications
in a variety of contexts.
Subsequently Koll\'{a}r proved that the singularities of
the theta divisors of a PPAV are mild in the sense that the pair 
$(A,\Theta )$ is log canonical. This implies for example that
$\Sigma _k(\Theta ):=\{ x\in A|{\rm mult }_x (
\Theta )\geq k\}$ is of codimension at least $k$ in $A$.
Ein and Lazarsfeld \cite{EL}
prove that if $\Theta $ is irreducible then $\Theta $ is normal and 
has only rational singularities.
In particular they show that if the codimension is exactly equal to $k$,
then $(A, \Theta)$ splits as a $k$-fold product of PPAVs.
In \cite{EL}, they also generalize the result of Koll{\'{a}}r
to $\QQ$-divisors. They prove that if $D$ is a divisor in the linear
series $|m\Theta |$ then the pair $(A,(1/m)D)$ is log canonical.
In \cite{H}, the first
result of Ein and Lazarsfeld was extended to the case of $\QQ$-divisors.
Let $(A,\Theta )$ be a PPAV, and
$D$ a divisor as above, such that $\lfloor (1/m)D\rfloor =0$. 
Then the pair $(A,(1/m)D )$ is log terminal.
These results immediately generalize to arbitrary polarizations.
Let $(X,L)$ be a polarization of type $(d_1,...,d_g)$, (as usual assume
that $d_i|d_{i+1}$), then there exists an isogeny $X\lra X_1$ 
of degree $d_g$, such that
$L$ is the pull back of a principal polarization $L_1$ on $X_1$.
If $D$ is a divisor in the linear series $|mL|$, then 
by pulling back by the corresponding isogeny $X_1\lra X$,
one deduces for example that the pair $(A,(1/md_g)D)$ is log canonical.

\vvv
\hhh In this paper, exploiting the theory of Fourier transforms,
we prove that a coherent sheaf which is generically of rank $1$ on an abelian
variety and is "cohomologically" a
principal polarization, must in fact be a line bundle and hence
a principal polarization.
A surprising consequence is that if $(X,L)$ is
a general polarized abelian variety
of type $(1,...,1,2)$, and $D$ is a divisor in the linear series
$|mL|$, then the pair $(A, (1/m)D)$ is log canonical.
When $(X,D)$ is an abelian surface of type $(1,2)$, and $D$ is irreducible, 
then the pair $(X,D)$ is log terminal. It is not clear whether
the last statement generalizes to $\QQ$-divisors.
Using the same technique, we are also able to give a cohomological 
criterion for a subvariety of an abelian variety to be a principal
polarization.
\vvv
\noindent \hskip 1cm {\bf Acknowledgment.}\hskip.3cm I would like to thank
R. Lazarsfeld, J. Koll\'{a}r, A. Bertram and A. Chen 
for the stimulating
conversations.
\vvv
\begin{center}{\bf 0. Notation and conventions}\end{center}
\vvv
 
$f^*D$          pull-back

%$D_{red}$        reduced divisor associated to $D$
%
%$\Omega ^i _X<D>$       sheaf of holomorphic $i$-forms with logarithmic poles along $D_{red}$
%
$|D|$    linear series associated to the divisor $D$

$Bs|D|$  the base locus of the linear series associated to the divisor $D$

$\equiv $ numerical equivalence

$\omega _X = \Omega ^n_X $  canonical sheaf of $X$ 

$\Omega ^i _X $
sheaf of holomorphic $i$-forms on $X$ 

$K_X$  linear equivalence class of a canonical divisor on $X$

$\omega _{X/Y}:=\omega _X \otimes (f^*\omega _Y)^*$
relative canonical sheaf of a morphism $f:X\lra Y$

$\kappa (X)$ Kodaira dimension of $X$

$k(x)$ the sheaf $\OO _X/m_{x}$

${\rm ALB} (Z)$ albanese variety of $Z$

${\rm alb}_Z: Z\lra {\rm ALB} (Z)$ albanese morphism

$\bar {Z}:={\rm alb}_Z (Z)$ albanese image of $Z$

$T^0$ connected component containing the origin of the subgroup $T$

\vvv
\noindent \hskip 1cm Unless otherwise stated $X$, 
$Y$ will denote smooth complex projective
varieties. 
If $D$ is a $\QQ$-divisor we will denote by $\lfloor D \rfloor $
and by $\lceil D \rceil$ the round down and the round up of $D$ respectively.
%If  $Z$ is a normal variety, by abuse of notation we will denote by
%$h^0(Z,\Omega ^i _Z)$, the integer $h^0(\tilde {Z},\Omega ^i _{\tilde {Z}})$
%for an appropriate smooth birational model $\tilde {Z}$ of $Z$. 

\vvv
\begin{center}{\bf 1. Preliminaries}\end{center}
\vvv
\noindent \hskip 1cm 
Let $X$ be a smooth complex projective variety and let $D$ be a 
$\mathbb Q$-divisor on $X$. We will say that $f:Y\lra X$ is a log-resolution 
of the pair $(X,D)$, if $f$ is a proper birational morphism such that 
$f^{-1}D \cup \{ \ exceptional\ set\ of \ f\ \}$
is a divisor with normal crossing support. Given a 
log-resolution of the pair $(X,D)$, we can define the multiplier ideal sheaf
associated to the divisor $D$,
$$\II (D):=f_* (\OO _Y (K_{Y/X}-\lfloor f^* D \rfloor ))\ .$$
The definition is independent of the choice of log-resolution.
Multiplier ideal sheaves may be defined in much grater generality. 
Their properties have been extensively studied 
eg. [N], [Sk], [De] and [E]. Under the above 
assumptions, we will say that the pair  $(X,D)$  is log canonical 
(respectively  log terminal) if the 
multiplier ideal sheaf associated to the divisor 
$(1-\epsilon )D$ is trivial, for $0<\epsilon < 1$ (respectively 
$0\leq \epsilon < 1$).
\vvv
\noindent \hskip 1cm 
One of the main tools that we will use is the theory of deformation
of cohomology groups developed by Green and Lazarsfeld \cite{GL1},
\cite{GL2}. For the
convenience of the reader, we include a brief summary of the 
main results. We wish to study the loci
$$V^i:= \{ P \in \Pic (X)\ |\ 
h^i(\omega _X\otimes P )\neq 0\}\ .$$
The geometry of these loci is governed by the following theorem
\cite{GL1}, \cite{GL2}, \cite{EL}:
\vvv
\noindent \hskip1cm {\bf Theorem 1.1 (Generic vanishing).}  {\it  
\vskip.2cm 

\noindent (i) $\Pic (X)\supset V^0(\omega _X  )\supset V^1(\omega _X )
\supset \ ...\ \supset V^n(\omega _X )=\{\OO _X \}$
\vskip.2cm 

\noindent (ii) Any irreducible component of $V_i$ is a 
translate of a torus and is
of codimension at least $i-(\dim (X) -\dim\ \alb _X(X))$ in $\Pic (X)$
\vskip.2cm 

\noindent (iii) Given $P$ a general point of an irreducible component $T$ of
$V^i$. Let $\phi \in H^0(X,\Omega ^1 _X)$.
Suppose $\bar {\phi } \in H^1(X, \OO _X)\cong T_{P }\Pic (X)$ is 
not tangent to $T$. Then the sequence
$$H^{0}(X,\Omega ^{n-i-1}_X\ot P) \stackrel{\wedge \phi} {\longrightarrow} 
H^{0}(X,\Omega ^{n-i}_X\ot P)  \stackrel {\wedge \phi } {\longrightarrow } 
H^{0}(X,\Omega ^{n-i+1}_X\ot P) $$
is exact. If $\bar {\phi }$
is tangent to $T$, then the 
maps in the above sequence vanish.}
\vvv
\noindent \hskip1cm {\bf Remark 1.2.}\hskip.3cm 
It is possible to generalize the
above definitions and results, to the more general case
of a morphism
$$\nu :X\lra A$$ 
from $X$ to an abelian variety $A$ see \cite {EL}.
One must then consider the loci
$$V^i(\omega _X,A):=\{P\in \Pic (A)|h^i(\omega _X\ot \nu ^*P)\neq 0\}.$$
The generic vanishing theorem still holds, with $\Pic (A) $ instead of
$\Pic (X)$ and $\dim (\nu (X))$ instead of $\dim ( \alb _X(X))$.
\vvv
\begin{center}{\bf 2. A characterization of principal polarizations}
\end{center} 

\vvv
\hhh In this section we will assume the notation and results of \cite {M}.
In particular $X$ will always be  an abelian variety and 
$\hat {X}$ will denote the corresponding dual
abelian variety. For any point $y\in \hat {X}$ let
$P_y$ denote the associated topological trivial line bundle.
Let $\PP$ be the normalized Poincar\'{e} bundle on
$X\times \hat {X}$. One may define a functor $\sh $ of $\OO_X$-modules
into the category of $\OO_{\hat {X}}$-modules by
$$\sh (M)=\pi _{\hat {X},*}(\PP \ot \pi ^*_XM).$$
The derived functor $R\sh$ of $\sh$ then induces an equivalence of 
categories between the two derived categories $D(X) $ and $D(\hat {X})$
\cite {M} theorem 2.2:

\vvv 
\noindent \hskip1cm {\bf Theorem 2.1 (Mukai).}{\it \hskip.3cm
There are isomorphisms of functors:
$$R\ss \circ R\sh \cong (-1_X)^*[-g]$$
and
$$R\sh \circ R\ss \cong (-1_{\hat {X}})^*[-g],$$
where $[-g]$ denotes "shift the complex $g$ places to the right".}
\vvv 
\hhh The index theorem (I.T.) is said to hold for a coherent sheaf $\FF$ on $X$
if there exists an integer $i(\FF )$ such that for all $j
\not \hskip-.15cm {=}i(\FF )$,
$H^j(X,\FF \ot P)=0$ for all $P\in \Pic (X)$.
The weak index theorem (W.I.T.) holds for a coherent sheaf $\FF$ if 
there exists an integer which we again denote by $i(\FF )$ such that for all $j
\not \hskip-.15cm {=}i(\FF )$, $R^j\sh
(\FF )=0$. It is easily seen that the I.T. implies the W.I.T.
We will denote the coherent 
sheaf $R^{i(\FF )}\sh (\FF )$ on $\hat {X} $ by $\hat {\FF }$.
One of the main themes of \cite{M} is that information on the cohomology
groups $H^i(X,\FF \ot P)$ may be interpreted as information
on the coherent sheaves $\FF$, $R^i \sh (\FF )$ and $\hat {\FF }$ 
(if $\FF$ satisfies
the W.I.T.). In this spirit we prove:

\vvv
\noindent \hskip1cm {\bf Proposition 2.2.}{\it \hskip.3cm
Let $X$ be an abelian variety, $\FF$ a coherent sheaf 
generically of of rank $1$ (eg. $\FF \cong L\ot \II$
where $L$ is a line bundle and $\II \subset \OO _X$ is an ideal sheaf).
For all $P\in \Pic (X)$, suppose that 
$h^0(X,\FF \ot P)=1$, and 
$h^i(X,\FF \ot P)=0$ for all $1\leq i\leq n=\rm {dim} (X)$.
Then $\FF$ is a line bundle with $h^0(X,\FF )=1$, hence a principal 
polarization.}
\vvv

\noindent {\it Proof.} \hskip.3cm We will use the
theory of Fourier functors developed by Mukai in \cite {M}. $\FF$ 
satisfies the I.T. and hence the W.I.T. and $i(\FF )=0$. 
We will need the following:

\vvv
\noindent \hskip1cm {\bf Claim 2.3.} {\it \hskip.3cm
$\hat {\FF }$ is a line bundle and $h^{n}(\hat{X},\hat {\FF })=1$.}
\vvv

\noindent {\it Proof of claim 2.3.}  \hskip.3cm
We will use the isomorphism of \cite {M} proposition 2.7
$$Ext^i_{\OO _X}({k(x)},\FF )\cong H^i(\hat{X},\hat {\FF }\ot P_{-x})$$
to compute $h^i(\hat{X},\hat {\FF })$.
Let $K^{p , q }$ be the bicomplex 
${\Cal {C}}^p[X,{\Cal {H}}om(k(x),\II ^q(\FF ))]$ where 
$\II ^{q }(\FF )$ is an injective resolution of the sheaf
$\FF $ (see \cite{G} II 7.3).
Since $\FF$ is generically of rank $1$, a local computation shows that
for generic $x\in X$, 
$$\EE xt ^{n}_{\OO _X}(k(x),\FF )\cong k(x),$$ 
and that for all 
$0\leq i \leq \dim (X)-1$
$$\EE xt ^{i}_{\OO _X}(k(x),\FF )\cong 0.$$ 
The corresponding spectral sequence therefore degenerates at the 
$E_2$ term and we have $'E^{p, q}_2=H^p(X,
\EE xt ^{q}_{\OO _X}(k(x),\FF ))$ and hence $'E^{p, q}_2\cong \CC $
for $p=0$ and $q=\dim (X)$, and $'E^{p,q}_2\cong 0$ otherwise.
Since the $E_{\infty }$ term gives a filtration of 
$Ext _{\OO _X}^{i }({k(x)},\FF )$, it follows that 
for $0\leq i\leq \dim (X)-1$ we have 
$Ext^i_{\OO (X)}({k(x)},\FF )\cong 0$.
On the other hand, for $i=\dim (X)$ 
$$Ext^{n}_{\OO (X)}({k(x)},\FF )\cong H^0(X,{k(x)} ) \cong \CC .$$
 The claim 
now follows since for a generic point $x\in X$ we have
$h^{n}(\hat{X},\hat {\FF }\ot P_{-x})=1$.
\hskip.5cm QED

\vvv
\noindent \hskip1cm Consequently $\hat {\FF }\cong \OO _{\hat {X}} (-\Theta )$ 
for an appropriate theta
divisor $\Theta$. $\hat {\FF }$ also satisfies the I.T.,
and in fact $h^i(\hat{X},\hat {\FF }\ot P_{-x})=0$ for all $0\leq i 
\leq \dim (X)-1$ and hence $\FF \cong (-1_X)^*\hat {\hskip -.07cm \hat {\FF}}$
is also a line bundle with only one section 
(see also \cite {M} proposition 3.11).
QED

\vvv
\noindent \hskip1cm {\bf Remark 2.4.} \hskip.3cm
The same proof also shows that if  $\FF$ is a coherent sheaf 
generically of rank $r$ such 
that for all $P\in \Pic (X)$, $h^0(X,\FF \ot P)=1$ and
$h^i(X,\FF \ot P)=0$ 
for all $1\leq i\leq n=\dim (X)$.
Then $\hat {\FF }$ is dual to an ample line bundle with
$r$ sections, and 
$\FF$ is a locally free sheaf. We will say that a sheaf 
satisfying the above properties is a cohomological
principal polarization. Our proposition states that a cohomological
principal polarization
of generic rank $1$
is a principal polarization.

\vvv
\noindent \hskip1cm {\bf Remark 2.5.} \hskip.3cm
The hypothesis on the cohomological groups may not be weakened.
Consider an irreducible theta divisor $\Theta \subset X$. Let
$f:\tilde {\Theta }\lra \Theta$ be an appropriate smooth birational model
Since $\Theta $ has only log-terminal singularities,
$h^0(f_*(\omega _{\tilde {\Theta }}\ot P))$
is $1$ if $P\not \hskip-.15cm {=}\OO_X$ and is $n$ if
$P=\OO _X$. Similarly for $1\leq i\leq n-1$,
$h^i(f_*(\omega _{\tilde {\Theta }}\ot P))$ is $0$ if 
$P\not \hskip-.15cm {=}\OO_X$ and is $\binom {n}{i+1}$ if
$P=\OO _X$. It follows that the sheaf $f_*(\omega 
_{\tilde {\Theta }})\oplus \OO_X$
is coherent, generically of rank $1$ and satisfies the 
conditions for being a cohomological principal polarization, for all $P
\not \hskip-.15cm {=}\OO_X$.
Analogously one could also consider the sheaf $k(0)\oplus \OO _X$.
It is not clear however whether one can find a similar example
for a sheaf of the form $L\ot \II$
where $L$ is a line bundle and $\II \subset \OO _X$ is an ideal sheaf.

\vvv
\noindent \hskip1cm {\bf Remark 2.6.} \hskip.3cm
One might expect that an analogue of proposition
2.2 might hold for more general polarizations. This is however not 
the case. Consider an abelian surface with polarization $L$ of type $(1,3)$.
It is easy to show that for any point $x$, and any topologically
trivial line bundle $P$, $h^i(L\ot \II _x \ot P)=0$ for all $i\geq 1$.
Of course it follows that the sheaf $L\ot \II _x$ satisfies the W.I.T., and in
general has all the cohomological properties of a polarization of type $(1,2)$.
However the sheaf $L\ot \II _x$ is not a line bundle.

\vvv
\begin{center}{\bf 3. Sub-varieties of abelian varieties}\end{center}
\vvv
\noindent \hskip1cm
Let X be an abelian variety of dimension $n$. 
We will say that a reduced irreducible divisor
$Z$ is a cohomological theta divisor if for any desingularization
$\nu :\tilde {Z}\lra Z$

i) $h^i(\omega _{\tilde {Z}})=\binom {n}{i+1}$ for all $0\leq i\leq n-1$

ii) $h^i(\omega _{\tilde {Z}}\ot \nu ^*P)=0$ for all $1\leq i\leq n-1$ and all
$\OO_X \not \hskip-.15cm {=}P\in \Pic (X)$.
\vvv
\noindent \hskip1cm
By the Hodge and Serre dualities, condition i) is equivalent to 
$h^0(\Omega ^j_{\tilde {Z}})=\binom {n}{j}$ for all $0\leq j\leq n-1$, and
condition ii) is equivalent to 
$h^0(\Omega ^j_{\tilde {Z}}\ot \nu ^* P)=0$ for all $0\leq j\leq n-2$
and all
$\OO_X \not \hskip-.15cm {=}P\in \Pic (X)$.
The above definition, immediately generalizes to the case of a map
$\nu : Z\lra X$ from a smooth $n-1$ dimensional variety $Z$ to 
a $n$ dimensional abelian variety $X$.
The following lemma will be useful:

\vvv
\noindent \hskip1cm {\bf Lemma 3.1.}{\it \hskip.3cm
Let $Z$ be a reduced irreducible divisor
of a $n$ dimensional abelian variety $X$. 
Let ${\tilde {Z}}$ be a smooth birational model of $Z$.
Then $Z$ is of general type if and only if $|\chi (\OO_{\tilde{Z}})|>0$.}

\vvv
\noindent {\it Proof.}\hskip.3cm 
If $Z$ is of general type, then $|\chi (\OO_{\tilde{Z}})|>0$ by \cite {EL}
theorem 3 and theorem 3.3. If $Z$ is not of general type, then 
one may consider the 
abelian subvariety $B:=\{ x\in X|x+Z\subset Z\}^{0}$.
By \cite{U}, we have 

i) $Z\lra Z/B$ is an \'{e}tale fiber bundle with fiber $B$,

ii) $Z\lra Z/B$ is birational to the Iitaka fibering of $Z$, and

iii) ${Z/B}$ is of general type.

\noindent We may assume that $\tilde {Z}\lra \widetilde {Z/B}$
is the Iitaka fibration. Restricting to a general fiber, since 
$(\omega _{\tilde{Z}} \ot P)_{|B}\cong P_{|B}$, it follows that 
$|\omega _{\tilde{Z}} 
\ot P|$ is empty unless 
$P_{|B}\cong \OO _B$. So $H^0(\tilde{Z},
\omega _{\tilde{Z}} \ot P)=0$
for general $P\in \Pic (X)$, and hence $|\chi (\OO_{\tilde{Z}})|=0$.\proofend

\vvv
\noindent \hskip1cm {\bf Lemma 3.2.}{\it \hskip.3cm
Let $\nu :Z\lra X$ be a morphism from a smooth variety
to an $n$ dimensional abelian variety.
Assume that $\nu$ restricted to $Z$ is generically finite onto 
a divisor $\bar {Z}=\nu (Z)$. Then

i) If $\bar {Z}$ is of general type, then $|\chi (\OO _Z)|
>0$. 

ii) If $|\chi (\OO _Z)|>0$, then $Z$ is of general type. 

iii) If $h^1(\omega _Z\ot \nu ^*P)=0$
for all but finitely many $P\in \Pic (X)$, then $\bar {Z}$
is of general type if and only if $Z$ is of general type.}

\vvv
\noindent {\it Proof.}\hskip.3cm
i) If $\chi (\omega _Z)=0$, then $h^0(\omega _Z\ot \nu ^*P)=0$ 
for general $P\in \Pic (X)$, so
$h^0(\omega _{\bar {Z}}\ot P)$ also vanishes for generic $P\in \Pic (X)$. 
Hence 
$\chi (\omega _{\bar {Z}})=0$, i.e. ${\bar {Z}}$ is not of general type.

\vvv
\hhh ii) 
If $Z$ is not of general type, then $\bar {Z}$ is not of general type.
Let $\beta$ be the generic fiber of the Iitaka fibration
$Z\lra Y$, then $\beta$ is \'{e}tale onto its image $B\subset X$.
In fact $\kappa (\beta )=\kappa (B)=0$, and hence $\beta$ and $B$ are 
abelian varieties of positive dimension. Now $(\omega _Z\ot \nu ^*P)
_{|{\beta }}=
\nu ^*P_{|\beta }$, and since $\beta \lra B$ 
is an \'{e}tale map of abelian varieties, 
we see that $h^0(\nu ^*P_{|\beta })=0$ for generic $P\in \Pic (X)$. 
So $h^0(\omega _Z\ot \nu ^* P)$ also vanishes for generic $P\in \Pic (X)$
and by the generic vanishing theorem $\chi (\omega _Z)=0$.

\vvv
\hhh iii) Assume now that $Z$ is of general type, 
and $h^1(\omega _Z\ot \nu ^*P)=0$ for all
but finitely many $P\in \Pic (X)$. If $\bar {Z}$ is not of general type, then 
$\chi (\OO _{\bar {Z}})=0$. Let $\tilde {Z}$ be a desingularization
of $\bar {Z}$. Since $\bar {Z}$ is a divisor of $X$, 
by \cite {EL} proposition 2.2, there  
exists a positive dimensional subgroup $T$ of $\Pic (X)$ such that for 
all $P\in T$, $h^0(\omega _{\tilde {Z}}\ot P)\geq 1$. 
Since $\chi (\OO _{\tilde {Z}})=0$, then also $h^0(\Omega ^{n-2} 
_{\tilde {Z}}\ot P
)=h^1(\omega _{\tilde {Z}}\ot P)\geq 1$ for $P\in T$.
It follows that also $h^0(\Omega ^{n-2} _{ {Z}}\ot \nu ^*P
)\geq 1$ for $P\in T$ which is a contradiction. \proofend

\vvv
\noindent \hskip1cm {\bf Theorem 3.3.}{\it \hskip.3cm
Let $\Theta $ be a reduced irreducible
divisor on an abelian variety $X$. Then $\Theta $
is a principal polarization if and only if 
${\Theta}$ is a cohomological theta divisor.}

\vvv
\noindent {\it Proof.}\hskip.3cm
We may assume that $(\tilde {X}, \tilde {\Theta})$ is a log resolution
of the pair $(X, \Theta )$. One has the exact sequence of sheaves:
$$0\lra \omega  _{\tilde {X}}\lra \omega  _{\tilde {X}} \ot 
\OO _{\tilde {X}}(\tilde {\Theta})
\lra \omega _{\tilde {\Theta}}\lra 0.$$
Pushing forward this sequence yields \cite {EL} the exact sequence of sheaves
on $X$:
$$0\lra \OO _{ {X}}\lra   
\OO _{{X}}({\Theta})\ot \II (\Theta )
\lra f_*\omega _{\tilde {\Theta}}\lra 0.$$
If $\Theta$ is an irreducible principal polarization, then 
by \cite{EL} the pair $(X, \Theta )$ is log terminal, so the multiplier ideal
sheaf $\II (\Theta )$ is trivial. The cohomology groups of 
$\omega _{\tilde {\Theta}}$ may now be easily computed from the second 
exact sequence, and we conclude that 
${\Theta}$ is a cohomological theta divisor.

\vvv
\hhh Assume now that $\Theta $ is a cohomological theta divisor.
The goal here is to show that the sheaf $\OO _{{X}}({\Theta})\ot \II (\Theta )$
is a cohomological principal polarization and then to apply proposition 2.2.
Of course this will follow from the equivalent statement for 
the sheaf $\omega _{\tilde {X}}\ot \OO _{\tilde {X}}(\tilde {\Theta})$.
For all $\OO_X \not \hskip-.15cm {=}P\in \Pic (X)$ 
and $0\leq i\leq n$, we have 
$h^i(\OO _X \ot P)=0$, so $H^i(X, \OO _{{X}}({\Theta})\ot \II (\Theta )\ot P)
\cong H^i(X,f_*\omega _{\tilde {\Theta}}\ot P)$. 
We therefore need only consider
the case $P=\OO _X$.
There is an injection $H^0(X,\Omega ^i_X)\lra H^0({\tilde {\Theta}},
\Omega ^i _{{\tilde {\Theta}}})$. Using the Hodge and Serre dualities, one
may translate this into the isomorphisms
$H^{n-i-1}({\tilde {\Theta}},\omega _{\tilde {\Theta}})\cong H^{n-i}(
\tilde {X},\omega _{\tilde {X}})$.
Consequently $h^0(\omega _{\tilde {X}}\ot 
\OO _{\tilde {X}}(\tilde {\Theta}))=1$
and $h^i(\omega _{\tilde {X}}\ot \OO _{\tilde {X}}(\tilde {\Theta}))=0$
for all $1\leq i \leq n$.\proofend

\vvv
\noindent \hskip1cm {\bf Corollary 3.4.}  {\it
Let $Z$ be a smooth variety mapping
finitely onto a divisor $\bar {Z}$ of an $n$ dimensional abelian variety $X$.
Assume that $\bar {Z}$ generates $X$, and that
$Z$ satisfies cohomological properties analogous
to those of a cohomological theta divisor. 
Then $\bar {Z}$ is a theta divisor.}

\vvv
\noindent \proof Let $\bar {Z}$ be the image of $Z$ by the map $\nu :Z\lra X$. 
By assumption $h^i(\omega _Z)=\binom {n}{i+1}$ for all
$0\leq i\leq n-1$ and $h^i(\omega _Z \ot \nu ^* P)=0$ for all 
$1\leq i\leq n-1$ and $\OO_X \not
 \hskip-.15cm {=}P\in \Pic (X)$.
The induced map
$H^0(X,\Omega ^1_X)\lra H^0(Z, \Omega ^1_Z)$ is an isomorphism.
Let $\tilde {Z}$ be a smooth birational model of $\bar {Z}$. We may assume that
$Z\lra \bar {Z}$ factors through $\tilde {Z}$. By lemma 3.2, $\bar {Z}$
is of general type. By \cite{KV} theorem 1, $h^i(\omega _{\tilde {Z}})\geq 
\binom {n}{i+1}$ for all $0\leq i\leq n-1$.
Clearly 
$$h^i(\omega _Z\ot \nu ^*P)=h^0(\Omega ^{n-i-1}_Z\ot \nu ^*P)\geq 
h^0(\Omega ^{n-i-1}_{\tilde {Z}}\ot P)=h^i(\omega _{\tilde {Z}}\ot P)$$ 
for all
$0\leq i\leq n-1$ and $P\in \Pic (X)$. 
It follows that ${ {Z}}$ is also a cohomological theta divisor.
Theorem
3.3 now proves the corollary.
\proofend

\vvv \noindent \hskip1cm If $Z$ is a reduced and 
irreducible divisor of an abelian variety,
then the cohomological theta divisor condition may be substantially weakened.

\vvv
\noindent \hskip1cm {\bf Corollary 3.5.}  {\it 
An irreducible reduced divisor $Z$ of an $n$ dimensional
abelian variety $X$ is a principal polarization if
and only if one of the following equivalent conditions holds

i) $Z$ is a cohomological theta divisor

ii) $Z$ is of general type, $h^0(\omega _Z )=n$, $h^1(\omega _Z \ot P)=0$ 
for all 
$\OO_X \not \hskip-.15cm {=}P\in \Pic (X)$

iii) $h^0(\omega _Z )=n$, $h^0(\omega _Z \ot P)=1$ for all 
$\OO_X \not \hskip-.15cm {=}P\in \Pic (X)$}

\vvv
\noindent \proof
The equivalence of condition $i)$ is proposition 2. Condition $i)$, clearly
implies the other two conditions (see lemma 3.1). 

\vvv
\hhh Assume that condition $ii)$ is satisfied.
Let ${\tilde {Z}}$ be a smooth birational model of $Z$.
By \cite {KV} theorem 1, it follows that $h^0(\Omega ^i _{\tilde {Z}})=
\binom {n}{i}$ for all $0\leq i\leq n-1$, and $|\chi (\OO _{\tilde {Z}})|=1$.
By \cite {EL} lemma 1.8, for any $P\in \Pic (X)$, the condition
$h^1(\omega _{\tilde{Z}} \ot P)=0$ implies 
$h^i(\omega _{\tilde {Z}} \ot P)=0$ for
 all $1\leq i\leq n-1$. Since $|\chi (\OO _{\tilde {Z}})|=1$, we
conclude that $h^0(\omega _{\tilde {Z}} \ot P)=1$ for all 
$\OO_X \not \hskip-.15cm {=}P\in \Pic (X)$. Thus
condition $ii)$ implies that $Z$ is a cohomological theta divisor.

\vvv
\hhh Assume now that condition $iii)$ 
is satisfied. Since $h^0(\omega _{\tilde {Z}}\ot P)
>0$ for generic $P\in \Pic (X)$, then $|\chi (\OO_{\tilde {Z}})|>0$
and by lemma 3.1, $Z$ is of general type.
We will now argue that condition ii) must also hold.
Suppose therefore that there exists some irreducible component $T$ of
$V^1:=\{ P\in \Pic (X)\ s.t.\ H^1({\tilde{Z}},\omega _{\tilde{Z}}\ot 
P)\not \hskip-.15cm {=}0\}$. By \cite {EL} or \cite {GL1}, \cite {GL2},
$T$ is a subtourus of $\Pic (X)$ of codimension at least $1$. Let 
$P\not \hskip-.15cm {=}\OO _Z$ be a general
point in $T$, let $\phi \in  H^0(X, \Omega ^1 _X)$ be any holomorphic $1$-form
which is not the pullback of a $1$-form on $S=T^*$. Consider the following 
complex of vector spaces:
$$H^0(X,\Omega ^{n-i-2}_{\tilde {Z}}\ot P)\stackrel 
{\wedge \phi }{\lra }H^0(X,\Omega ^{n-1-i}_{\tilde {Z}}\ot P)
\stackrel {\wedge \phi }{\lra }H^0(X,\Omega ^{n-i}_{\tilde {Z}}\ot P).$$
By \cite {EL}, this is exact for all $i\geq 1$
(notice that any component of $V^i$, $i\geq 2$
through the general point $P$ is also contained in $T$ and hence may
be assumed to be either equal to $T$ or empty.
We may of course also assume that $P$ is a
general point of $V^i$, $i\geq 2$). Since $\chi$ is zero on exact sequences,
and $|\chi (\OO _X)|=1$ it follows that $H^0 (X,\Omega ^{n-2}_{\tilde {Z}}
\ot P) \stackrel {\wedge \phi }{\lra } 0$.
We may choose local 
coordinates $\{ z_1,...,z_{n-1}\}$ on $Z$, 
centered at a general point $z\in Z$. 
Choose $\phi _i\in H^0 (X,\Omega ^{1}_X)$ 
such that $dz_1,...,dz_{n-1}$ correspond to the restriction 
of $\phi _i$ to $T_z(Z)^*$. 
Let $\eta \in H^0(X,\Omega ^1 _X)$ be a non zero holomorphic 1 form such that
$\eta _{|T_z(Z)^*}=0$.
Since $Z$ is of general type, $Z$ is not vertical with 
respect to 
$\pi :X\lra S$, and for general $z\in Z$, 
$\eta$ is not in $\pi ^* (H^0(S,\Omega ^1
_S))$. We may further assume that the forms $\phi _i$ are not pulled back
from $S$. (If this is not the case, 
then consider instead $\phi _i+\epsilon _i \eta$.) 
By a local computation, it follows that given any non
zero $\omega \in H^0 
(X,\Omega ^{n-2}_{\tilde {Z}})$, 
there exists a $\phi _i$ such that $\omega \wedge
\phi_i$ is not zero. Since this is a contradiction, 
$h^1(Z,\omega _{\tilde {Z}}\ot P)=0$ for all $\OO_X
\not \hskip-.15cm {=}P\in \Pic (X)$.\hskip.5cm \proofend

\vvv
\begin{center}{\bf 4. Polarizations of type $(1,...,1,2)$}\end{center}

\vvv
\noindent \hskip1cm {\bf Theorem 4.1.}  {\it 
Let $(A,Z )$ be an $n$ dimensional
abelian variety with a polarization of type
$(1,...,1,2)$ (i.e. $Z$ is the class of an ample line bundle with two 
sections), and let
$D$ be any divisor in the linear series $|mZ |$. Then either 

i) $(X,L)$ splits as the product
of a PPAV and an elliptic curve, or 

ii) $(A, ( {1}/{m})D)$ is log canonical.}

\vvv
\noindent {\it Proof.}\hskip.3cm Consider the exact sequence associated to the
multiplier ideal sheaf $\II := \II (( 1-(\epsilon / m) )D)$ for an 
appropriate rational number $0<\epsilon  <<1$:

$$0\lra L\ot \II \lra L \lra L\ot \OO _X/\II \lra 0.$$
Since $h^i(L\ot \II \ot P)=0$ for all $1\leq i\leq n$ and all $P\in \Pic (X)$,
it follows that the above sequence is exact on global sections, and that
$h^0(L\ot \II \ot P)=\chi (L\ot \II )$ for all $P\in \Pic (X)$.
Since $L$ is of type $(1,...,1,2)$, we must have $0\leq \chi (L\ot \II )
\leq 2$. If $\chi (L\ot \II )=0$ then 
$h^i(L\ot \II \ot P)=0$ for all $0\leq i\leq n$ and all $P\in \Pic (X)$,
and hence by a result of Mukai $L\ot \II =0$ (as sheaves) which is 
impossible.
If $\chi (L\ot \II )=1$, then $L\ot \II$ is a cohomological
principal polarization and hence 
$L\ot \II$ is a line bundle with only one section. This is again a 
impossible, unless $\II \cong \OO_X (-Z ')$ where
$Z'$ is a divisor of type $(0,...,0,1)$. In this case
there is a map $X\lra E$ of $X$ onto some elliptic curve $E$
such that $Z'$ is the pull back of a principal polarization
on $E$. Let $L':=L\ot \OO_X (-Z')$. Then $(X,L')$ is a PPAV 
that splits as the product of
$(E,\OO _E (p))$ and the PPAV $(\bar {X},\bar {L})$. 
Finally if $\chi (L\ot \II )=2$, 
all sections of $L$ vanish along any translate of 
the cosupport of $\II$. This is only possible if the cosupport of $\II$
is empty, i.e. if $\II =\OO _X$.\proofend
\vvv
\noindent \hskip1cm {\bf Remark 4.2.}\hskip.3cm
This statement is optimal in the sense that one may consider
$X=E_1\times E_2$, the product of two elliptic curves, endowed with the 
polarization of type $(1,2)$ given by $Z =p_1^*(p)+p_2^*(q_1+q_2)$
where $p$ and $q_i$ are points of $E_1$ and of $E_2$ respectively.
If the $q_i$ are distinct points, then $(X,Z )$ is log canonical but
not log terminal. If $q_1=q_2$, then $(X,Z)$ is not log canonical.
One may ask whether for $D$ a reduced and irreducible
divisor of type $(1,...,1,2)$, the pair $(X, D)$ is log terminal.
While this is not clear in higher dimensions, it does hold if 
$X$ is an abelian surface. 

\vvv
\noindent \hskip1cm {\bf Theorem 4.3.}{\it \hskip.3cm
Let $(X,D )$ be a polarized abelian surface with $D$ an ample reduced
irreducible divisor such that $h^0(\OO_A(D))=2$. Then  
the pair $(X, D)$ is log terminal.}

\vvv
\noindent \proof As in the proof of proposition 2, one may consider
an appropriate smooth birational model $f:\tilde {D}\lra D$ and the 
induced exact sequence of sheaves
$$0\lra \OO _{ {X}}\lra   
\OO _{{X}}(D)\ot \II (D )
\lra f_*\omega _{\tilde {D}}\lra 0.$$

\noindent Assume that $(X,D)$ is not log terminal, then the cosupport of
$\II :=\II (D )$ is not empty and has codimension two in $X$. 
Pick a general divisor $D'\in |D|$
which is distinct from $D$. These two divisors will intersect along 
a finite set of points. It follows that $h^0(\OO _X(D)\ot \II \ot P)=2$
only if the appropriate translate of the cosupport of $\II$ 
is contained in this finite set of points. (By \cite{LB} 10.1.2, $|D|$
has exactly four base points, corresponding to $D'\cap D''$ the
intersection of two general members of $|D|$.) This means that 
$h^0(\OO _X(D)\ot \II \ot P)\leq 1$ for all but finitely many $P\in \Pic (X)$.
For all topologically trivial line bundles $P \not \hskip -.2cm {=}
\OO_X$ there is an isomorphism $H^0(\OO _{{X}}(D)\ot \II (D )\ot P)\cong
H^0(f_*\omega _{\tilde {D}}\ot P)$. 
Since $D$ is ample and hence of general type,
$h^0(f_*\omega _{\tilde {D}}\ot P)>0$ for all $P\in \Pic (X)$. So 
$h^0(f_*\omega _{\tilde {D}}\ot P)=h^0(\OO _{{X}}(D)\ot \II (D )\ot P) =1$ 
for all but 
finitely many $P\in \Pic (X)$. 
Hence $g({\tilde {D}})=2$ and
$h^0(f_*\omega _{\tilde {D}}\ot P)=h^0(\OO _{{X}}(D)\ot 
\II (D )\ot P) =1$ for all
$\OO _X$
$\not \hskip -.15 cm {=}P\in \Pic (X)$.  By corollary 3.5 the proposition
now follows. 
%Alternatively, for a direct proof, one may proceed as follows.
%The ideal sheaf $\II $
%is cosupported in dimension $0$. If the cosupport $\II $ consists of at least 
%two points, then $2g(\tilde {D})-2=D^2-\sum m_i (m_i-1)\leq 0$ and then
%$D$ is either a rational or an elliptic curve. In either case this is a 
%contradiction. So $\II $ is cosupported at a unique point.
%Furthermore $D\cap D'$ must also consist of a single point,
%otherwise $h^0(\OO _{{X}}(D)\ot \II (D )\ot P)$ would be equal to $2$
%for some nontrivial $P$.
%In which case $\OO _X(D)$ would have a unique base point at which
%every divisor in the linear series $|\OO _X(D)|$ vanishes with multiplicity
%exactly $2$, and any two divisors intersect exclusively along this point.
%This is a contradiction, eg. \cite{LB} 10.1.3.
\hskip.5cm QED
\vvv
\noindent \hskip1cm {\bf Question 4.4.}{\hskip.3cm Let $(X,D)$ be
a general polarized abelian of type $(1,...,1,2)$. Assume that $D$
is irreducible. Is $(X,D)$ a log terminal pair? Similarly let $(X,D)$ be
a general polarized abelian of type $(1,...,1,d)$, with $d=3,4$.
Assume that $D$
is irreducible. Is $(X,D)$ a log canonical pair?}
\vvv
%%%%%%%%%%%%%%%%%%%%%%%%%%%%%%%%%%%%%%%
% References
%%%%%%%%%%%%%%%%%%%%%%%%%%%%%%%%%%%%%%%

\begin{center} --------------------------------

\vskip.4cm
Department of Mathematics, 
University of Utah, Salt Lake City, Utah 84112-0090, USA. 

\vskip.3cm
email: hacon@math.utah.edu\end{center}

\begin{thebibliography}{99}

\bibitem {AM} {\it A. Andreotti, A. Mayer},
{{On period relations for abelian integrals on algebraic curves}},
Ann. Scuola Normale Superiore Pisa {\bf 21} (1967) 189-238


\bibitem {EL} {\it L. Ein, R. Lazarsfeld},
{{Singularities of theta divisors, and the birational geometry of
irregular varieties.}}, Jour. Amer. Math. Soc. {\bf 10} (1997), 243-258 


\bibitem{G} {\it R. Godement},
{{Topologie alg\'{e}brique et théorie des faisceaux}}, Publications de 
l'Institut de Math\'{e}matique de l'Universit\'{e} de Strasbourg, XIII. 
Actualit\'{e}s Scientifiques et
Industrielles, No. 1252. Hermann, Paris, 1973

\bibitem{GL1} {\it M. Green, R. Lazarsfeld}, {{Deformation theory, 
generic vanishing theorems, and some conjectures
of Enriques, Catanese and Beauville}}, Inventiones Math. {\bf 90} 
(1987), 389-407 

\bibitem{GL2}{\it  M. Green, R. Lazarsfeld},
{{Higher obstructions to deforming cohomology groups of line bundles}},
J. Amer. Math. Soc. {\bf 4} (1991), 87-103

\bibitem {H} {\it C. D. Hacon},
{{Divisors on principally polarized abelian varieties}},
Compositio Mathematica (to appear)

\bibitem {Ke} {\it G. Kempf},
{{On the geometry of a theorem of Riemann}},  
Ann. of Math. {\bf 98} (1973), 178-185

\bibitem{KV} {\it Y. Kawamata, E. Viehweg}, {{On a characterization of
abelian varieties in the classification theory of algebraic varieties}},
Compositio Math.  {\bf 41} (1980), 355-360

\bibitem {LB} {\it H. Lange C. Birkenhake},
{{Complex Abelian Varieties}}, 
Grundlehren der mathematischen Wissenschaften, Springer Verlag


\bibitem{M} {\it S. Mukai},
{{Duality between $D(X)$ and $D(\hat {X})$, with application to Picard sheaves }},
Nagoya math. J. {\bf 81} (1981), 153-175 

\bibitem{U}{\it  K. Ueno},
{{Classification Theory of Algebraic
Varieties and Compact Complex Spaces}},
Springer Verlag {\bf LNM 439 } (1975) 
\end{thebibliography}
\end{document}